\documentclass[10pt]{amsart}

\usepackage[text={400pt,650pt},centering]{geometry}

\usepackage{color}
\usepackage{esint,amssymb}
\usepackage{graphicx}
\usepackage{MnSymbol}
\usepackage{mathtools}
\usepackage[colorlinks=true, pdfstartview=FitV, linkcolor=blue, citecolor=blue, urlcolor=blue,pagebackref=false]{hyperref}
\usepackage{microtype}

\definecolor{darkgreen}{rgb}{0,0.5,0}
\definecolor{darkblue}{rgb}{0,0,0.7}
\definecolor{darkred}{rgb}{0.9,0.1,0.1}

\newtheorem{theorem}{Theorem}

\theoremstyle{definition}

\numberwithin{equation}{section}
\numberwithin{theorem}{section}

\newcommand{\R}{\mathbb{R}}

\newcommand{\LL}{\mathbb{L}}
\newcommand{\e}{\varepsilon}

\newcommand{\I}{\mathcal{I}}
\newcommand{\Hh}{\mathcal{H}}

\newcommand{\Msym}{M_{\mathrm{sym}}^{2\times 2}}
\newcommand{\Md}{M^{2\times 2}}

\newcommand{\per}{\mathrm{per}}

\begin{document}

\title[Isotropy and loss of strong ellipticity through homogenization]{Isotropy prohibits the loss of strong ellipticity through homogenization in linear elasticity}

\begin{abstract}
Since the seminal contribution of Geymonat, M\"uller, and Triantafyllidis, it is known that strong ellipticity is not necessarily conserved by homogenization in linear elasticity.
This phenomenon is typically related to microscopic buckling of the composite material. The present contribution is concerned with the interplay between isotropy and strong ellipticity in the framework of periodic homogenization in linear elasticity.
Mixtures of two isotropic phases may indeed lead to loss of strong ellipticity when arranged in a laminate manner.
We show that if a matrix/inclusion type mixture of isotropic phases produces macroscopic isotropy, then strong ellipticity cannot be lost.

\medskip

\noindent {\sc R\'esum\'e.} Nous savons depuis l'article fondateur de Geymonat, M\"uller et Triantafyllidis qu'en \'elasticit\'e lin\'eaire l'homog\'en\'eisation p\'eriodique ne conserve pas n\'ecessairement l'ellipticit\'e forte.
Ce ph\'enom\`ene est li\'e au flambage microscopique des composites. Notre contribution consiste ˆ examiner le r\^ole de l'isotropie dans ce type de pathologie. Le m\'elange de deux phases isotropes peut en effet conduire \`a cette perte si l'arrangement est celui d'un lamin\'e. Nous montrons qu'en revanche, si un arrangement de type matrice/inclusion produit un tenseur homog\'en\'eis\'e isotrope, alors la forte ellipticit\'e est conserv\'ee.
\end{abstract}

\author[G. A. Francfort]{Gilles A. Francfort}
\address[G. A. Francfort]{L.A.G.A.  (UMR CNRS 7539),  Universit\'e Paris-Nord, Villetaneuse, France}
\email{gilles.francfort@univ-paris13.fr}

\author[A. Gloria]{Antoine Gloria}
\address[A. Gloria]{Universit\'e Libre de Bruxelles (ULB), Brussels, Belgium and Inria, Lille, France}
\email{agloria@ulb.ac.be}

\subjclass[2010]{}
\date{\today}

\maketitle

%%%%%%%%%%%%%%%%%%%%%%%%%%%%%%%%%%%%%%%%%%%%%%
%%%%%%%%%%%%%%%%%%%%%%%%%%%%%%%%%%%%%%%%%%%%%%

\section{Introduction}

This contribution is restricted to the two-dimensional case, although most of the results that are being used remain true in any dimension.

In all that follows $\mathfrak T$ stands for the unit 2-torus. Consider a $\mathfrak T$-periodic heterogeneous linear elastic material characterized by its elasticity tensor field $\LL$, a $\mathfrak T$-periodic  symmetric endomorphism on $\Msym$, the set of $2\times2$-symmetric matrices.
Assume that $\LL$ is an $L^\infty$, pointwise very strongly elliptic map, that is that, for some $\lambda>0$, 
\begin{equation}\label{vse}M \cdot \LL(x) M \ge \lambda |M|^2 \text{ for  all $M\in \Msym$, and a.e. $x\in \mathfrak T$}.
\end{equation}
Let $D$ be an open,   bounded domain of $\R^d$, and $u_0 \in H^1(D;\R^2)$.
Then it is simple to establish, in view of  classical homogenization results, that the integral functional  $$H^1_0(D)\ni u \mapsto \I_\e(u):=\int_D \nabla (u+u_0) \cdot \LL(\frac{x}{\e})\nabla (u+u_0)dx$$ $\Gamma$-converges for the weak topology of $H^1_0(D)$ 
to the homogenized integral functional $$H^1_0(D)\ni u \mapsto \I_*(u):=\int_D \nabla (u+u_0) \cdot \LL_* \nabla (u+u_0)dx,$$ where $\LL_*$ is a constant elasticity tensor that is very strongly elliptic with constant $\lambda$. It is classically given by
\begin{equation}\label{def_hom}
M\cdot \LL_*M  =\min\left\{\int_{\mathfrak T}(M+\nabla v)\cdot\LL(x) (M+\nabla v)\ dx; \ v\in \Hh^1({\mathfrak T};\R^2)\right\}.
\end{equation}

If instead of pointwise very strong ellipticity, we only assume pointwise strong ellipticity, that is that, for some $\lambda>0$,  $$M \cdot \LL(x) M \ge \lambda |M|^2 \text{ for  all  symmetrized {\it rank one}  $M=a\otimes b,\ a,b\in\R^2$, and a.e. $x\in {\mathfrak T}$,}$$ the story is different. 

In an inspirational work \cite{GMT-93}, G. Geymonat, S. M\"uller, and N. Triantafyllidis introduced  two measures of coercivity:
\begin{eqnarray*}
\Lambda\;\;\;\;&=&\inf\left\{\frac{\int_{\R^d}{\nabla u}\cdot \LL\nabla u\ dx}{\int_{\R^d}|\nabla u|^2\ dx}; u\in C^\infty_0(\R^2,\R^2)\right\},
\\[2mm]
\Lambda_{\per}&=& \inf\left\{ \frac{{\int_{\mathfrak T} {\nabla v} \cdot \LL \nabla v\ dx}}{{\int_{\mathfrak T} |\nabla v |^2\ dx}};\,v\in \Hh^1({\mathfrak T};\R^2)\right\}.
\end{eqnarray*}

When $\Lambda>0$, then homogenization occurs as in the classical setting of \eqref{vse}, while when $\Lambda<0$, $\I_\e$ is not bounded from below and there is no homogenization. Their focus was on the case when $\Lambda=0$ and $\Lambda_{per}>0$. There, they showed that there is still homogenization towards $\I_*$ with associated $\LL_*$ given by \eqref{def_hom}. However, $\LL_*$  may be strongly elliptic (that is, non-degenerate on rank-one matrices) or simply non-negative on such matrices, 
 but not strongly elliptic because there exists 
$a, b\in \R^2$ such that $a \otimes b\cdot \LL_* a\otimes b=0$).

The third phenomenon is referred to as \emph{loss of strong ellipticity by homogenization}.
To avoid confusion we will say that a fourth-order tensor $\LL$ is strongly elliptic if $M\cdot \LL M\ge 0$ for all rank-one matrices, and that it is \emph{strictly} strongly elliptic if in addition there exists $\lambda>0$ such that this inequality can be strengthened to $M\cdot \LL M\ge \lambda |M|^2$.

There is only one single example \cite{Gutierrez-98}  for which one can prove that strong ellipticity is lost by homogenization.
The associated composite material has a laminate structure made of two isotropic phases (a strong phase and a weak phase). 
Loss of strong ellipticity occurs when the strong phase buckles in compression (it is related to the failure of the cell-formula for nonlinear composites, cf. \cite{Muller-87,GMT-93}). This has been rigorously proved in \cite{Briane-Francfort-15}.

Buckling is by nature a very one-dimensional phenomenon. It is mechanically  unlikely  that any material could lose strong ellipticity in \emph{every} rank-one direction.
This simple-minded observation suggests that assuming the isotropy of $\LL_*$ may prevent loss of strong ellipticity by homogenization. 
The aim of the present contribution is precisely to mathematically corroborate the mechanical intuition.

Let us quickly describe our main result.
 The fact that $\LL_*$ is isotropic allows one to focus on the Lam\'e coefficients of $\LL_*$. If they held true, the Hashin-Shtrikman (HS) bounds would prevent loss of strong ellipticity \emph{a priori}. Whereas the HS bound on the bulk modulus  does hold  true, we do not know whether even the elementary harmonic lower bound for the Poisson's ratio similarly  holds true. The standard proof for \emph{very} strongly elliptic materials  proceeds by duality and cannot be used in our setting of strongly elliptic materials since the energy density is not necessarily pointwise non-negative.
Instead, we argue through a comparison argument  which does not use duality.

%\begin{remark}
%As pointed out to us by A. Braides and M. Briane, elementary abstract results make it possible to prove the $\Gamma$-convergence result of \tref{alaGMT} under the sole assumption
%of $\Lambda\ge 0$. We have kept the condition $\Lambda_6>0$ for its crucial role in the rest of this article.
%\end{remark}
%

\section{The result} \label{sec:ishom}

Let  $\lambda_1,\mu_1$ and $\lambda_2,\mu_2$ be the
Lam\'e coefficients of isotropic stiffness tensors $\LL^1$ and $\LL^2$. In other words,
$$(\LL^i)_{pqrs}=\lambda_i \delta_{pq}\delta_{rs}+\mu_i(\delta_{pr}\delta_{qs}+\delta_{ps}\delta_{qr}),\; i,p,q,r,s\in\{1,2\}.$$
We assume that
\begin{equation}
\label{eq:cond}
0<\mu_1=-(\lambda_2+\mu_2)=:-K_2<\mu_2,\; K_1:=\lambda_1+\mu_1>0.
\end{equation}
In particular $\lambda_i+2\mu_i\ge 0,\ i=1,2,$ so that both phases are strictly strongly elliptic but phase 2 is not very strongly elliptic because $K_2<0$.

We then define
$$\LL(x)= \chi(x)\LL^1+(1-\chi(x))\LL^2,$$
where $\chi$ is the  characteristic function the inclusion (phase 1),  an open subset of the torus ${\mathfrak T}$ with Lipschitz boundary;  assume further that
\begin{equation}\label{scon}
\{\chi=0\}:=\{x\in {\mathfrak T};\ \chi(x)=0\} \text{ is  connected in } {\mathfrak T}.
\end{equation}

The following result is a generalization  of \cite[Theorem~2.9, Case 2]{Briane-Francfort-15} because, in contrast with that result,  it does not impose any restriction on the geometry of each phase besides representing the worst inclusion/matrix type microstructure, that is that for which the matrix (here phase 2) does not satisfy very strong ellipticity. 
\begin{theorem}\label{ploom}
Under assumptions \eqref{eq:cond}, \eqref{scon},  $\Lambda \ge 0$
and $\Lambda_\per>0$.
\end{theorem}
As shown in \cite[Proposition~3.4]{Briane-Francfort-15} the laminate configuration (in the  periodic setting) results in a loss of ellipticity for $\LL_*$.
The theorem below shows that the isotropy of $\LL_*$ \emph{rules out  any loss of ellipticity}.
\begin{theorem}\label{t.HS}
Under assumptions \eqref{eq:cond}, \eqref{scon}  assume further that  
$\LL_*$ is isotropic with  bulk modulus $K_*$ and shear modulus  $\mu_*$. 
Then $K_*+\mu_*>0$, that is, $\LL_*$ is strictly strongly elliptic.
\end{theorem}
We expect that Theorem~\ref{t.HS} also holds in the stationary ergodic setting (for which statistical isotropy is a mild requirement that yields isotropy of $\LL_*$).
The proof we display below fails however to cover this setting due to the use of Korn's theorem on the (compact!) torus.

%%%%%%%%%%%%%%%%%%%%%%%%%%%%%%%%%%%%%%%%%%%%%%
%%%%%%%%%%%%%%%%%%%%%%%%%%%%%%%%%%%%%%%%%%%%%%

\section{Proofs}

\subsection{Proof of Theorem \ref{ploom}}
\label{3.1}

\hskip0.5cm

\medskip

\noindent \textit{Step~1.} $\mathbf{\Lambda\ge 0}$.
We decompose $\LL$ as $\LL-\underline \LL+\underline \LL$ where $\underline \LL$ is the isotropic stiffness tensor with Lam\'e constants $\underline \lambda,\underline \mu$
defined as follows: $\underline\mu:=\mu_1$ and $\underline \lambda:=\inf_x\{\lambda(x)+\mu(x)\}-\mu_1= -2\mu_1$.
On the one hand, so defined, $\underline \LL$ is clearly strongly elliptic since $\underline \mu=\mu_1>0$ and $\underline \lambda+2\underline \mu= 0$.

On the other hand, 
$
\mu-\underline \mu\ge 0,\;(\lambda+\mu)-(\underline \lambda+\underline \mu) \ge0,
$
so that 
\begin{equation}\label{non-neg}\LL-\underline \LL \text{ is pointwise non-negative as a quadratic form.}
\end{equation}
This yields $\LL\ge \underline \LL$ pointwise, which implies that $\Lambda\ge 0$ since $\underline \LL$ is  strongly elliptic, so the corresponding $\Lambda$ is non-negative.

\medskip

\noindent \textit{Step~2.} $\mathbf{\Lambda_\per>0}$.
We argue by contradiction and assume that $\Lambda_\per=0$. Consider $\nabla v_n$ a minimizing sequence of periodic fields  with $\int_{\mathfrak T} |\nabla v_n|^2\ dx=1$ such that
\begin{equation}\label{e.ass-to-zero}
\lim_{n\to \infty}\int_{\mathfrak T} \nabla v_n\cdot \LL\nabla v_n\ dx\,=\,0.
\end{equation}
We now prove that 
(i) $\nabla v_n \rightharpoonup 0$ weakly in $L^2({\mathfrak T};\R^{2\times2})$, then that
(ii) $\nabla v_n$ is strongly convergent in  $L^2({\mathfrak T};\R^{2\times2})$.
The combination of (i) and (ii) then yields $\lim_{n\to \infty} {\int_{\mathfrak T} |\nabla v_n|^2\ dx}=0$, whence the contradiction.

The proof of (i) exploits the structure of the problem.
In the spirit of the proof of  \cite[Theorem~2.9]{Briane-Francfort-15} we add a null Lagrangian $4\mu_1 \det \nabla v_n$ (which satisfies $\int_{\mathfrak T} \det \nabla v_n \ dx\equiv 0$) to the energy, so that the assumption turns into
\begin{equation}\label{e.conv-zero-energy}
\lim_{n\to \infty}\int_{\mathfrak T} \{\nabla v_n\cdot \LL\nabla v_n+4\mu_1 \det \nabla v_n\}\,dx=\,0.
\end{equation}
On the one hand, since $\int_{\mathfrak T} |\nabla v_n|^2\ dx= 1$, we may assume that (along a subsequence) there exists a periodic field $\nabla v$   in $L^2({\mathfrak T};\R^{2\times2})$ such that 
$\nabla v_n\rightharpoonup \nabla v$ weakly in $L^2({\mathfrak T};\R^{2\times2})$. Since $\Lambda \ge 0$, the map $\nabla u \mapsto \int_{\mathfrak T} \nabla u\cdot \LL\nabla u\ dx$ is weakly lower-semicontinuous, which implies that
$$
{\int_{\mathfrak T} \nabla v\cdot \LL\nabla v\ dx}=0.
$$
As in Step~2 of the proof of  \cite[Theorem~2.9]{Briane-Francfort-15}, we have,
pointwise for all $u=(u^1,u^2)\in \Hh^1({\mathfrak T},\R^2)$, 
$$
\nabla u\cdot \LL\nabla u+4\mu_1 \det \nabla u\,=\, 
P(\frac{\partial u^1}{\partial y_1},\frac{\partial u^2}{\partial y_2})
+R(\frac{\partial u^1}{\partial y_2},\frac{\partial u^2}{\partial y_1}),
$$
where $P$ and $R$ are  quadratic forms that satisfy,
for some $\alpha>0,$
\begin{eqnarray}
P(a,b)&\geq & \alpha (a+b)^2 \chi+\alpha(a-b)^2 (1-\chi),\label{eq:strong1}\\
R(a,b)&\geq & \alpha (a-b)^2\chi+\alpha(a^2+b^2)(1-\chi).\label{eq:strong2}
\end{eqnarray}
In particular, since these quadratic forms are non-negative, they vanish
almost everywhere  at $\nabla v$.
On the set $\{\chi=1\}$, this yields
\begin{eqnarray}
\frac{\partial v^1}{\partial y_1}+\frac{\partial v^2}{\partial y_2}&=&0,\label{eq:strong-correction1}\\
\frac{\partial v^1}{\partial y_2}-\frac{\partial v^2}{\partial y_1}&=&0\label{eq:strong-correction2},
\end{eqnarray}
while, on the set $\{\chi=0\}$, 
\begin{eqnarray}
\frac{\partial v^1}{\partial y_1}&=&\frac{\partial v^2}{\partial y_2},\label{eq:strong-correction3}\\
\frac{\partial v^1}{\partial y_2}&=&\frac{\partial v^2}{\partial y_1}=0.\label{eq:strong-correction4}
\end{eqnarray}
From \eqref{eq:strong-correction2} and \eqref{eq:strong-correction4}, we deduce that there exists a potential $\psi\in \Hh^2({\mathfrak T},\R)$ such that
$v=(v_1,v_2)=\nabla \psi$.
We start with proving additional properties on $\psi$ in the matrix, that is on the set $\{\chi=0\}$. By assumption, this set is  connected, so that
from \eqref{eq:strong-correction4} we deduce that $\psi(y)=\psi_1(y_1)+\psi_2(y_2)$ for some $\psi_1,\psi_2 \in H^2({\mathfrak T};\R)$ on that set.
From \eqref{eq:strong-correction3} we then learn that $\psi_1(y_1)=ay_1^2+by_1+c$ and that $\psi_2(y_2)=ay_2^2+dy_2+e$ for some
$a,b,c,d,e\in \R$.
We continue with the properties of $\psi$ in the inclusion, that is on the set where $\{\chi=1\}$.
On the one hand, taking the derivative of \eqref{eq:strong-correction1} w.~r.~t. $y_1$ and of \eqref{eq:strong-correction2}  w.~r.~t. $y_2$, and using the Schwarz' commutation rule, we obtain
that $-\Delta v_1\,=\,0$. On the other hand, the formula for $\psi$ on the set where $\{\chi=0\}$ completes this equation with the boundary data $v_1(y)=2ay_1+b$.
By uniqueness of the solution of this boundary-value problem, we then conclude that $v_1(y)=2ay_1+b$ on ${\mathfrak T}$. Likewise, $v_2(y)=2ay_2+d$.
In turn the condition $\int_{\mathfrak T} \nabla v\ dx=0$ due to periodicity  implies that $\nabla v \equiv 0$  as claimed.

\medskip

We turn now to the proof of (ii) and shall argue that if \eqref{e.ass-to-zero} holds, then $(1-\chi)\nabla v_n$ converges strongly in $L^2({\mathfrak T};\R^{2\times2})$ to zero.
Integrating \eqref{eq:strong2} over the unit torus ${\mathfrak T}$ yields in view of \eqref{e.ass-to-zero}
\begin{equation}\label{int1}
\frac{\partial v^1_n}{\partial y_2}\to 0,\;  \frac{\partial v^2_n}{\partial y_1}\to 0, \text{ strongly in } 
L^2(\{\chi=0\};\R^{2\times 2}).\end{equation}

\bigskip
It remains to prove that  ${\partial v^1_n}/{\partial y_1}$ and
${\partial v^2_n}/{\partial y_2}$ converge strongly to zero in $L^2(\{\chi=0\}; \R^{2\times2})$ as well.
By symmetry it is enough to treat the first term  ${\partial v^1_n}/{\partial y_1}$. To this aim we follow the beginning of the argument of Step~2 in the proof of  \cite[Theorem~2.9]{Briane-Francfort-15}. We get
$$\frac{\partial v^1_n}{\partial y_1}-\frac{\partial v^2_n}{\partial y_2}\to 0, \text{ in } L^2(\{\chi=0\}; \R^{2\times2}).$$

 Consequently, with \eqref{int1},
$$
\begin{array}{l}
%\item $\frac{\partial v^1_n}{\partial y_1}$ converges  in $L^2(\Omega,L^2(Z_2))$,
\displaystyle\frac{\partial}{\partial y_2}\frac{\partial v^1_n}{\partial y_1}=\frac{\partial}{\partial y_1}\frac{\partial v^1_n}{\partial y_2}\to 0 \text{ in }H^{-1}(\{\chi=0\}),\\[3mm]
\displaystyle \frac{\partial}{\partial y_1}\frac{\partial v^1_n}{\partial y_1}=\frac{\partial}{\partial y_1}\frac{\partial v^2_n}{\partial y_2}+ r_n (\text{ with } r_n\to 0 \text{ in }H^{-1}(\{\chi=0\}))=\frac{\partial}{\partial y_2}\frac{\partial v^2_n}{\partial y_1}+ r_n\to 0\text{ in }H^{-1}(\{\chi=0\}).
\end{array}
$$
As in the proof of  \cite[Theorem~2.9]{Briane-Francfort-15}, application of  Korn's theorem \cite{Necas} then yields
$$
 \frac{\partial v^1_n}{\partial y_1}\to 0, \text{ strongly in }L^2(\{\chi=0\};\R^2).
 $$

Likewise ${\partial v^2_n}/{\partial y_2}$ converges strongly to zero in $L^2(\{\chi=0\};\R^2)$, and, in view of \eqref{int1}$,\nabla v_n$ converges strongly to zero in $L^2(\{\chi=0\};\R^2)$, or, equivalently,
\begin{equation}\label{eq:strong3}
\int_{\mathfrak T} (1-\chi) |\nabla v_n|^2 \ dx\to 0.
\end{equation}

We are now in a position to conclude. Since $\det \nabla v_n$ is quadratic in $\nabla v_n$, \eqref{eq:strong3} yields $(1-\chi)\det \nabla v_n\to 0$ in $L^1({\mathfrak T})$ while $\int_{\mathfrak T}\det \nabla v_n\ dx= 0$, 
so that
\begin{equation}\label{int2}
\int_{\mathfrak T} \chi \det \nabla v_n\ dx\to 0.
\end{equation} 
From \eqref{eq:strong1}, \eqref{eq:strong2}, 
$$
\int_{\mathfrak T} \chi \left(  \left\{\frac{\partial v_n^1}{\partial y_1}+\frac{\partial v_n^2}{\partial y_2}\right)^2+\left(\frac{\partial v_n^1}{\partial y_2}-\frac{\partial v_n^2}{\partial y_1}\right)^2 \right\}\ dx \to 0.
$$
Thus, with \eqref{int2}, 
$$
\int_{\mathfrak T} \chi |\nabla v_n|^2\ dx=\int_{\mathfrak T} \chi \left(  \left\{\frac{\partial v_n^1}{\partial y_1}+\frac{\partial v_n^2}{\partial y_2}\right)^2+\left(\frac{\partial v_n^1}{\partial y_2}-\frac{\partial v_n^2}{\partial y_1}\right)^2 \right\}\ dx -2\int_{\mathfrak T} \chi \det \nabla v_n\ dx\to 0.
$$
Combined with \eqref{eq:strong3}, this yields the desired contradiction 
since ${\int_{\mathfrak T} |\nabla v_n|^2}\ dx=1$.

%%%%%%%%%%%%%%%%%%%%%%%%%%%%%%%%%%%%%%%%%%%%%%
%%%%%%%%%%%%%%%%%%%%%%%%%%%%%%%%%%%%%%%%%%%%%%

\subsection{Proof of Theorem~\ref{t.HS}}

We split the proof into two steps and prove $(K_*,\mu_*)\ge (-\mu_1,\mu_1)$ and $(K_*,\mu_*)\neq (-\mu_1,\mu_1)$ separately.

\medskip

\noindent \textit{Step~1.} 
Since $\Lambda\ge 0$ and $\Lambda_\per>0$, \eqref{def_hom} defines the homogenized elasticity tensor, so that
for all $M\in \Md$ there exists $v_M\in \Hh^1({\mathfrak T};\R^2)$ such that 
$$
M\cdot \LL_* M\,={\int_{\mathfrak T} (M+\nabla v_M)\cdot \LL (M+\nabla v_M)\ dx}.
$$
Consider $\underline{\LL}$ as in Subsection \ref{3.1}.
Since $\underline \LL$ is constant and  strongly elliptic  while $\int_{\mathfrak T} \nabla v_M\ dx=0,$ 
$$
\int_{\mathfrak T} (M+\nabla v_M)\cdot \underline \LL (M+\nabla v_M)\ dx=\, M \cdot \underline \LL  M + \underbrace{\int_{\mathfrak T} \nabla v_M \cdot \underline \LL \nabla v_M\ dx}_{\displaystyle \ge 0} \,\geq \,  M \cdot \underline \LL  M.
$$
%.
Hence, 
$$
M\cdot (\LL_*-\underline \LL) M \,\geq \, {\int_{\mathfrak T} (M+\nabla v_M)\cdot (\LL-\underline \LL) (M+\nabla v_M)\ dx}.
$$
Appealing to \eqref{non-neg}, we conclude that $M\cdot (\LL_*-\underline \LL) M\ge 0$. The isotropy assumption on $\LL_*$ then permits to conclude that   $(K_*,\mu_*)\ge (-\mu_1,\mu_1)$.

\medskip

\noindent \textit{Step~2.}  We argue by contradiction and assume that $(K_*,\mu_*)= (-\mu_1,\mu_1)$, so that $\LL_*=\underline\LL$.
In this case all rank-one matrices  $a\otimes a$ are such that 
\begin{equation}\label{int3} a\otimes a\cdot \LL_*a\otimes a=0.
\end{equation}
Then,
$$
0 = \int_{\mathfrak T} (a\otimes a+\nabla v_{a\otimes a}) \cdot \LL(a\otimes a+\nabla v_{a\otimes a})\ dx.
$$
Adding the null-Lagrangian $4\mu_1 \det (a\otimes a+\nabla \phi_{a\otimes a})\ dx$ and proceeding as in Step~2 of the proof of Theorem \ref{ploom}, 
we conclude that $\phi_{a\otimes a}\equiv cst.$
Hence the homogenization formula takes the form
$$
a\otimes a\cdot \LL_* a\otimes a \,=\, a\otimes a\cdot \left[\int_{\mathfrak T} \LL(x)\ dx\right]a\otimes a >0
$$
since the volume fraction of phase 1 is not $0$. This contradicts \eqref{int3} and concludes the proof.

\section*{Acknowledgements}
AG acknowledges financial support from the European Research Council under
the European Community's Seventh Framework Programme (FP7/2014-2019 Grant Agreement
QUANTHOM 335410).

\bibliographystyle{plain}
%\bibliography{strongellipticity}

\def\cprime{$'$}

\end{document}